\documentclass[11pt]{article}
\usepackage{amsfonts,amssymb,amsmath,amsthm}
\usepackage{graphicx}
\usepackage{bbm}
\usepackage{psfrag}
\usepackage{calc}
\usepackage{epic}
\usepackage[english]{babel}
\usepackage{hyperref}
\usepackage{float}
\usepackage{latexsym}
\usepackage{color}
\usepackage{cite}
\usepackage{subfigure,dsfont}

\newtheorem{theorem}{Theorem}[section]

\newtheorem{corollary}[theorem]{Corollary}
\newtheorem{remark}[theorem]{Remark}

\newcommand{\be}{ \begin{equation}}
\newcommand{\ee}{\end{equation}}
\newcommand{\ben}{ \begin{equation*}}
\newcommand{\een}{\end{equation*}}

\def\E{{\mathbb E}}
\def\I{{\mathbbm I}}
\def\P{{\mathbb P}}
\def\L{{\mathcal L}}

\newcommand{\IsDef}{\triangleq}			

\newcommand{\footremember}[2]{%
    \footnote{#2}
    \newcounter{#1}
    \setcounter{#1}{\value{footnote}}%
}

\begin{document}

\title{
The end time of SIS epidemics driven by random walks on edge-transitive graphs}

\author{%
  Daniel Figueiredo\footremember{UFRJ1}{Systems Engineering and Computer Science, Federal University of Rio de Janeiro (UFRJ), Rio de Janeiro, Brazil, daniel@land.ufrj.br}%
  \and
Giulio Iacobelli\footremember{UFRJ2}{Statistical Methods Department, Federal University of Rio de Janeiro (UFRJ), Rio de Janeiro, Brazil, giulio@im.ufrj.br}%
\and
  Seva Shneer\footremember{HW}{Heriot-Watt University, V.Shneer@hw.ac.uk}%
  }

\maketitle

\begin{abstract}
Network epidemics is a ubiquitous model that can represent different phenomena and finds applications in various domains. Among its various characteristics, a fundamental question concerns the time when an epidemic stops propagating. We investigate this characteristic on a SIS epidemic induced by agents that move according to independent continuous time random walks on a finite graph: Agents can either be infected (I) or susceptible (S), and infection occurs when two agents with different epidemic states meet in a node. After a random recovery time, an infected agent returns to state S and can be infected again. The End of Epidemic (EoE) denotes the first time where all agents are in state S, since after this moment no further infections can occur and the epidemic stops.
For the case of two agents on edge-transitive graphs, we characterize EoE as a function of the network structure by relating  the Laplace transform of EoE to the Laplace transform of the meeting time of two random walks. Interestingly, this analysis shows a separation between the effect of network structure and epidemic dynamics. We then study the asymptotic behavior of EoE (asymptotically in the size of the graph) under different parameter scalings, identifying regimes where EoE converges in distribution to a proper random variable or to infinity.  
We also highlight the impact of different graph structures on EoE, characterizing it under complete graphs, complete bipartite graphs, and rings.


\end{abstract}

\section{Introduction}
\label{sec:intro}

Network epidemic models are an ubiquitous and powerful abstraction that can represent different phenomena in various domains, such as physics, biology, and social sciences. The classic model assumes network nodes correspond to individuals and edges indicate the possibility of direct influence. In this model, nodes have an epidemic state that changes over time according to some function of the epidemic state of their respective neighbors. The most elementary epidemic state is represented by a single binary digit, and thus, every node is found in one of two possible states, frequently denoted by susceptible (S) and infected (I)~\cite{pastor2015epidemic,draief2010epidemics,barabasi2016network}.

A fundamental problem concerning network epidemics is understanding the final (or time average) epidemic state of the nodes as time unfolds. Intuitively, the network structure plays a key role, as illustrated by the celebrated work of Pastor-Satorras and Vespignani, showing that the epidemic threshold vanishes on networks where the degree distribution is heavy enough~\cite{pastor2001epidemic}. Indeed, the role of network structure on specific epidemic models has been broadly investigated and different dichotomies have been identified (e.g., very long versus very short epidemic duration)~\cite{draief2010epidemics,ganesh2005effect,moez2006epidemic,newman2002spread,pastor2015epidemic}. Moreover, recent efforts have focused on understanding the impact of dynamic network structure (i.e., edge set changes over time)~\cite{masuda2017introduction,tunc2013epidemics,valdano2015analytical}.

Another class of network epidemic models consider agents that move on the network. In this model, network nodes represent locations where agents can reside and edges indicate the possibility of direct movement between locations. In addition, the epidemic state is now associated with the agents (and not nodes) and  changes over time when agents with different states meet in a node. Note that this model embodies two different dynamics, namely agent  mobility and epidemic diffusion. While epidemic diffusion clearly depends on agent mobility, agent mobility may be independent of epidemic diffusion. Nevertheless, this coupled dynamics adds significant complexity, making a rigorous theoretical analysis much more challenging. Indeed, most theoretical results on this model are fairly recent when compared to the classic network epidemic model\cite{datta2004random,dickman2010activated,draief2011random,abdullah2011viral,benjamini2016epidemic,kesten2005spread,nagatani2018}. However, this model has been considered and analyzed through numerical simulations for at least 45 years~\cite{kelker1973}, since it also finds applications in various domains. 

Arguably the simplest agent mobility model are random walks, where agents choose neighbors uniformly at random and independently from one another. Indeed, this is the preferred choice in theoretical works that tackle this model. Moreover, the simplest kind of epidemic is the SI model, where every agent has a binary state (S or I) and can only transition from state S to I. For example, Draief and Ganesh~\cite{draief2011random} consider an SI epidemic with two random walks and characterize the infection probability over time as a function of the network, illustrating again the importance of the network structure on the epidemic. In a more recent work, Nagatani et al. analyze the SIS model (where agents alternate between S and I states) with many independent walkers (metapopulation model) on different networks to show that infection risk and epidemic threshold are a function of the network structure~\cite{nagatani2018}. 

In the SIS epidemic, an agent in state I (infected) returns to state S (susceptible) after some time, known as recovery time, and can become infected again. However, the epidemic stops when all agents are found in state S, as no agent can further become infected. Let the end of the epidemic (EoE) denote the first time instant where all agents are found in state S. Under mild conditions (finite graph, finite number of walkers, recovery time with finite moments), EoE is finite almost surely. However, its value strongly depends on model parameters and network structure. Thus, EoE is a crucial quantity of SIS dynamics as it reveals a fundamental property of the epidemic, namely, when it ends. This metrics has been investigated in different models, as discussed in Section~\ref{sec:related}.

The main contribution of this work is a characterization of EoE as a function of the network structure. We consider edge-transitive graphs and two independent random walks with exponentially distributed step time and recovery time. Under this assumption, we provide the exact Laplace transform for EoE as a function of the Laplace transform for meeting times (Theorem~\ref{thm:main}). Interestingly, the graph structure only influences the latter which does not depend on the epidemic dynamics. On an intuitive level, our main result separates the effect of the network structure from the epidemic dynamics. 

Our second contribution is the characterization of EoE on graph sequences of increasing size. In particular, we identify scaling regimes for which EoE converges to a distribution (with finite moments) or diverges to infinity (Theorems~\ref{thm:limit_general} and~\ref{thm:limit_general_2}). Interestingly, while on fixed graphs EoE is finite, graph sequences with a proper scaling allow the EoE to grow with the graph size. Moreover, the scaling regimes necessary for EoE to diverge strongly depend on the graph, again illustrating the importance of network structure. We illustrate this behavior by considering complete graphs, complete bipartite graphs, and rings.

The remainder of this paper is organized as follows. A  summary of  related work and results is presented in Section~\ref{sec:related}. Section~\ref{sec:notation} presents the notation and some preliminary definitions. The main result is given in Section~\ref{sec:main} providing the Laplace transform for the EoE for any edge-transitive graph. The main result is later applied in Section~\ref{sec:limit} to derive limit results and scaling regimes for different graphs, along with auxiliary theorems to characterize their behavior. Last, a brief discussion and outlook concerning this problem is presented in Section~\ref{sec:conclusion}.

\section{Related Work}
\label{sec:related}

In what follows the two different network epidemic models are presented more formally along with some of important results on the characterization of the epidemic.

\subsection{Epidemics on nodes}

In this class of network epidemics nodes correspond to individuals in a given population and edges encode the possible interactions among the population. The various epidemic states such as S (susceptible), I (infected), and R (recovered), are associated with the nodes, and change over time according to the epidemic state of neighboring nodes. In such models, epidemic dynamics is strongly driven by network structure with node degree playing a fundamental role. 

One of the most famous models in this class is a SIS epidemic model normally refereed to as {\em contact process}, defined as follows. Let $G=(V,E)$ denote an undirected graph with node and edge set $V$ and $E$, respectively, $\lambda$ and $\delta$ two parameters called the contagion rate and recovery rate, respectively (usually $\delta=1$). In this model the state of a node $i$ evolves as follows: 
 \begin{align*}
 \underset{\text{\em susceptible / infected} }{\xi_i(t)=\{0,1\}}  \qquad  \begin{cases}
 \xi_i:0\to 1 &  \text{ \em at rate $\lambda \sum_{(j,i)\in E}\xi_j$}\\
  \xi_i:1\to 0 & \text{ \em at rate $\delta$} 
 \end{cases}
 \end{align*}
The model was introduced on an infinite lattice by Harris~\cite{harris1974contact} 45 years ago, and it has been broadly explored. There are several surveys on this topic providing many details and generalizations of this classic model~\cite{durrett1995ten,liggett2013stochastic}. Note that besides the network structure, $\lambda$ plays a fundamental role. Intuitively, if $\lambda$ is much larger than $\delta$, infections occur much faster than recovery and the epidemic may spread very quickly. Below we provide a few important results concerning the survival of the epidemic in different scenarios. 

\medskip
\noindent
{\bf On infinite lattices:}
When $G = \mathbb{Z}^d$, there exists an {\em epidemic threshold} $\lambda_c:=\inf\{\lambda: \mathbb{P}(\xi(t) \neq 0\, \forall t>0)>0\}$ such that 
\begin{itemize}
\item 
$\lambda<\lambda_c$: epidemic dies out a.s., that is  $\mathbb{P}(\exists t_0: \xi(t)= 0\, \forall t>t_0)=1$,  $\forall \xi(0)$

\item
$\lambda > \lambda_c$:  epidemic survives with positive probability (at any node), that is 
 \[
 \mathbb{P}(\xi(t) \neq 0 ,\forall t>0)>0 \text{ and } \;\forall i \quad   \mathbb{P}(\forall T \exists t>T: \xi_i(t)=1 )>0\;,
 \]
\end{itemize}
for all initial conditions with infinitely many infected nodes.

\medskip
\noindent
{\bf Regular infinite trees}: 
In this case, two epidemic thresholds $\lambda_1<\lambda_2$ have been identified. For $\lambda<\lambda_1$ and for $\lambda>\lambda_2$ the behavior is identical to $\mathbb{Z}^d$. Moreover,
\begin{itemize}
\item
 $\lambda \in (\lambda_1, \lambda_2)$: epidemic survives with positive probability,  but every node recovers eventually a.s., that is 
\[ 
\mathbb{P}(\xi(t) \neq 0, \forall t>0)>0 \text{ and } \forall i \quad   \mathbb{P}(\exists T: \xi_i(t) =0 \forall t>T)=1 
\]
\end{itemize}
where the initial condition is one node called the root infected~\cite{pemantle1992contact}.

\medskip
\noindent
{\bf On arbitrary finite graphs:}
In this case it is known that the epidemic will eventually die out with probability one, independently of the network structure and infection rate. However, there is still a phase transition on the time that the epidemic ends, which can be either very early or very late. Let $G=([n],E)$ be an undirected finite and connected graph on $n$ vertices,  $A$ its adjacency matrix, and $\rho(A)$ its spectral radius. Let $\tau$ denote the time that the epidemic ends, defined as follows:
\[ 
 \tau= \inf\{t>0: \xi_i(t)=0\;, \forall i \in [n]\}.
 \]
Then the following holds:
\begin{itemize}
\item 
$\lambda \rho(A) < 1$  $\implies$ $\mathbb{E}(\tau)\leq \frac{\log n +1}{1-\lambda \rho}$,  $\forall \; \xi(0)$,

\item $\lambda \eta(G)>1$  $\implies$ $\exists C>0$ such that $\mathbb{E}(\tau)\geq e^{Cn}$, $\forall \xi(0)\neq 0$, where $\eta(G)= \inf_{S\subset [n]: |S|\leq \lfloor n/2\rfloor} \frac{E(S, S^c)}{|S|}$ is the {\em isoperimetric constant} 
and $E(S,S^c)$ denotes the number of edges between $S$ and $S^c$.
\end{itemize}
Thus, the expected time for the epidemic to end can be logarithmic or exponential in the size of the network, depending on its structure and infection rate~\cite{ganesh2005effect,moez2006epidemic}.

\subsection{Epidemics on agents}

In this class of network epidemic models, nodes correspond to locations and the edges encode the possibility of direct movement between locations. Different from the previous model, network nodes have no epidemic state. In contrast, the epidemic state is associated with agents that move around on the network. The epidemic state of agents can change when agents with different states meet in a network node. Moreover, it is often assumed that agents perform independent random walks in continuous time. Note that in this model network structure influences agent mobility which in turn influences the epidemic dynamics, adding an extra layer of complexity with respect to the previous model. 

\medskip
\noindent
{\bf SI epidemic on finite graphs with two agents:}
Consider two agents performing independent random walks $X_t$, $Y_t$, in continuous time according to rate transition matrix $Q$ (assumed to be  reversible). Start with one agent infected and the other susceptible and no recovery. Assume the susceptible agent becomes infected as a function of the time it has spent together with the infected agent, on any given node. 
Let $\mathbbm{I}(t)$ be an indicator for the susceptible agent to be infected by time $t$, and $p(t)$ its expectation. Then 
\begin{itemize}
\item 
$\mathbb{E}(\tau(t))=\sum_{i \in V} \pi_i^2 t$
\item
$p(t)\leq 1- e^{-\beta t \sum_{i \in V} \pi_i^2}$
\end{itemize}
where $\tau(t)$ is the {\em coincidence time} up to time $t$ (total time the two agents have spent together, up to time $t$), and $\pi$ is the invariant distribution for the random walk on $G$. Note that graphs with different degree distribution (i.e., regular graph versus power law distribution) will have very different scalings for this quantities~\cite{draief2010epidemics}. 

\medskip
\noindent
{\bf SI epidemic on infinite lattices:}
Let $G = \mathbb{Z}^d$ and consider two types of agents all performing independent continuous time random walks: $A$-particles ({\em susceptible}) step with rate $D_A$ and $B$-particles ({\em infected}) step with rate $D_B$. Assume there is no recovery and infection occurs immediately when an $A$-particle meets a $B$-particle in a node. 

Let $N_A(x,0_-)\sim Poi(\lambda)$ denote the number of $A$-particles at $x\in \mathbb{Z}^d$ at time $0_-$. Moreover, let $N_B$ denote a fixed number of $B$-particles placed in the lattice (not at random) at time $0$. Define the two sets:
\[
B(t):=\left\{ x \in \mathbb{Z}^d: \text{ a $B$-particle visits $x$ during $[0,t]$}\right\} + \left[-\frac{1}{2}, \frac{1}{2}\right]^d
\]
\[
\mathcal{C}(r):= [-r, r]^d
\]
Then, for any $D_A, D_B\geq 0$ there exists a  constant $C_1<\infty$ (independent  of $N_B$ and initial positions of $B$-particles) such that for  all sufficiently large $t$
\[
\mathbb{E}(\text{number of B-particles outside $\mathcal{C}(C_1t)$ at time t})\leq 2 N_B e^{-t}
\]
$\qquad \implies$ $B(t)\subset \mathcal{C}(2C_1 t)$ eventually almost surely~\cite{kesten2005spread}.


Moreover, if $D_A = D_B$ then there exists a constant $C_2>0$ such that for  each constant  $K>0$ and  for sufficiently large $t$
\[
\mathbb{P}(\mathcal{C}(C_2t)\not\subset B(t))\leq \frac{1}{t^K}
\]
$\qquad \implies$  $\mathcal{C}(C_2t)\subset B(t)$ eventually almost surely~\cite{kesten2005spread}.

These two theorems suggests that  a ``shape theorem'' may hold: $t^{-1}B(t)$ converges to a non-random set $B_0$ which implies that the growth rate of $B(t)$ is linear in $t$. Note that when $D_A=0$ which implies that susceptible particles do not move, the model degenerates to what is known as the {\bf frog model}~\cite{alves2002phase}. In this case there exists a full shape theorem, as follows: $\exists$ a non random set $B_0$ such that for all $0<\varepsilon<1$ 
{
\[
(1-\varepsilon)B_0 \subset \frac{B(t)}{t}\subset (1+\varepsilon)B_0\, ,\qquad  \text{for all $t$ large enough a.s.}
\]
This means that
$\mathbb{P}(\forall \varepsilon \exists t_0(\varepsilon): \forall t>t_0\,\; (1-\varepsilon)B_0 \subset \frac{B(t)}{t}\subset (1+\varepsilon)B_0) =1.$~\cite{alves2002phase}

\medskip
\noindent
{\bf SI epidemic on finite graphs:}
The frog model has also been recently studied on finite graphs~\cite{benjamini2018epidemic}. Let $G=(V,E)$, $N_A(x, 0_-)\sim Poi(\lambda)$, $N_B=1$ a single infected particle, placed in a given node. In this model, $A$-particles do not move (frog model) and $B$-particles perform independent discrete time random walks. However, $B$-particles have a lifetime $\tau$ (not random) after which they are removed from the system. The process stops when there are no more $B$-particles present. Let $R_\tau$ denote the set of nodes visited by $B$-particles (with lifetime $\tau$) before the process stops. Let $S(G) = \inf\{t>0: R_t = V\}$ denote the smallest lifetime required by $B$-particles to visit every node of $G$. A recent work has characterized the asymptotic behavior of $S(\mathbb{T}_d(n))$ as a function of the walking rate, where $\mathbb{T}_d(n)$ are regular trees with $n$ nodes~\cite{hermon2018frogs}.

\medskip
\noindent
{\bf SIS epidemic on infinite lattices:}
The SIS epidemic where agents are driven by the frog model has also been investigated on infinite graphs. In this model, a $B$-particle becomes an $A$-particle when an $A$-particle moves into its node (an S to I transition), while an $A$-particle becomes a $B$-particle after an exponential amount of time, with rate $\lambda$. However, only $A$-particles move according to independent random walks. The time at which the epidemic ends (i.e., all agents are found in state S) has been characterized, showing a phase transition between very short and very long~\cite{dickman2010activated,sidoravicius2017absorbing}. In particular, a phase transition on the density of the agents has recently been shown for infinite lattices of any dimension, $\mathbb{Z}^d$ for fixed $d>0$~\cite{sidoravicius2017absorbing}. 

\medskip
\noindent
{\bf SIS epidemic on finite graphs with two agents:}
This is the scenario tackled in this paper, focusing on edge-transitive graphs, that is described in the following sections. Note that this is the first rigorous work on characterizing the end time of SIS epidemics on finite graphs with mobile agents, to the best of our knowledge.

\section{Notation and preliminaries}
\label{sec:notation}

Let $G=(V,E)$ be a undirected, finite, connected graph with $|V|=n$ vertices, and let $d(i)$ denote the degree of vertex $i \in V$. Throughout the paper we assume that $G$ is edge-transitive, i.e., given any two edges $e_1, e_2\in E$, there exists an automorphism of $G$ that maps $e_1$ to $e_2$. 
Informally speaking, the edge transitivity assures that every edge ``sees'' the same graph structure. Notable examples are: complete graph ($K_n$), complete bipartite graph ($K_{n_1,n_2}$), cycle graph ($C_n$), star graph ($S_n$) and the hypercube ($Q_n$) with $d>0$ dimensions (where $n=2^d$).

Consider two agents moving on the graph $G$ according to independent continuous time random walks, denoted by $\{W_1(t)\}_{t\geq 0}$ and  $\{W_2(t)\}_{t\geq 0}$, with $W_k(t)\in V$, $k=1,2$,  for any time $t$. For each walker, the holding time in every vertex is exponentially distributed with parameter $\lambda$ ({\em walking rate}), independently from the other walker. Thus, an agent in vertex $i$ moves to vertex $j$ with rate $\lambda/d(i)$ if $\{i,j\} \in E$, and $0$ otherwise. 

We assume agents are either {\em susceptible} (S) or {\em infected} (I), denoted by $\{S_k(t)\}_{t\geq 0}$ with $S_k(t) \in \{S,I\}$ for $k=1,2$. We consider a SIS epidemic. When an infected agent meets the susceptible in a vertex, an infection immediately occurs. Note that this event takes place when the S agent walks into the vertex where I resides, or when the I agent walks into the vertex where S resides. Once infected, an agent recovers by transitioning to the S state after some time. The recovery time is assumed to be exponentially distributed with parameter $\gamma$ ({\em recovery rate}), and is
independent of other events. Once in the S state, the agent becomes infected again when it meets the I agent. 

The system dynamics can be fully described by the joint state of both agents, $(S_k(t), W_k(t))$ for $k=1,2$. We assume both agents are infected and located in the same vertex at time zero, and thus $S_1(0)=S_2(0)=I$ and $W_1(0)=W_2(0)=i$, for some $i\in V$. 

The following are important quantities related to this model:
\begin{itemize}
\item[-]
{\bf Meeting time} of the two walkers. Let $i,j \in V$ and define 
\begin{align}\label{eq:tau}
O^{i,j}:= \inf \{ t\geq  0: W_1(t) =W_2(t) \mid W_1(0)=i \text{ and } W_2(0)=j  \}\;.
\end{align}
\end{itemize}
A worst case polynomial upper bound (in $n$) for the meeting of two walkers in any graph is shown in~\cite{coppersmith1993collisions}. An upper bound for the expected meeting time for a fixed $G$, given in terms of the hitting time of a single walker, is also known~\cite{aldous2002reversible}. The Laplace transform of this meeting time has also been established in closed form for some specific graphs, including a scenario with more than two walkers~\cite{ohwa2015}.

A notion related to the meting time above is the time for two walkers to meet when they start at distance one from each other. This is a fundamental quantity for the analysis of our model, as we soon discuss.
\begin{itemize}
\item[-]
{\bf Meeting time from distance one} of the two walkers. In this case, we restrict $i,j \in V$ such that $(i, j) \in E$, and thus,
\begin{align}\label{eq:M}
M^{i,j}:= \inf \{ t\geq  0: W_1(t) =W_2(t) \mid W_1(0)=i \text{ and } W_2(0)=j \text{ and } (i, j) \in E \}\;.
\end{align}
\end{itemize}
Note that since $G$ is edge-transitive the distribution of $M^{i,j}$ does not depend on the specific edge $(i, j) \in E$, and all edges have the same distribution. Henceforth, we drop the indication of the edge in the notation and denote by $M$ the random variable with this distribution. 


The previous quantities depend only on the graph structure and walking rate $\lambda$, but not on the epidemic dynamics. In sharp contrast, the following quantity indicates the time that the epidemic ends which occurs when both walkers become susceptible. 
\begin{itemize}
\item[-] {\bf End of epidemic time} is defined as 
\begin{align}\label{eq:T}
T:= \inf \{t\geq 0: S_1(t)=S_2(t)=S\}\;,
\end{align}
for a given initial condition $W_k(0) = i$ and $S_k(0) = I$, for $k=1,2$ and $i \in V$.
\end{itemize}
Note that both agents will stay in the susceptible state ever after time $T$, moving on the graph but never becoming infected again. Moreover, note that $T$ depends on graph structure $G$, the walking rate $\lambda$ and the recovery rate $\gamma$, while $M$ does not depend on the recovery rate $\gamma$. In what follows we provide a characterization of $T$ as a function of $G$, $\lambda$, and $\gamma$.

\section{End of Epidemic Time}
\label{sec:main}

We now state our main theorem relating the End of Epidemic time $T$ to the meeting time from distance one $M$.
\begin{theorem} \label{thm:main}
Let $\L_T(s) = \E(e^{-sT})$ denote the Laplace transform of $T$ and $\L_{M}(s) = \E(e^{-sM})$ denote the Laplace transform of $M$. Then, for any $s > 0$,
\begin{align} \label{eq:Laplace2}
\L_T(s) &= \frac{2\gamma \left(\frac{1 - \L_M(s+\gamma)}{s+\gamma} -  \frac{1 - \L_M(s+2\gamma)}{s+2\gamma}\right)}{\frac{2 \lambda + s}{2\lambda} - 2 \L_M(s+\gamma) + \L_M(s+2 \gamma)}
\end{align}
\end{theorem}

{\bf Proof of Theorem~\ref{thm:main}.}

Denote by $R_k \sim Exp(\gamma)$ the time to recovery of individual  $k$ and define  $R_* \IsDef \min(R_1,R_2) \sim Exp(2 \gamma)$ the shortest time to recovery. Denote  by $J_k \sim Exp(\lambda)$ the time to the next jump of walker $k$ and define $J_* \IsDef \min(J_1,J_2) \sim Exp(2 \lambda)$ the shortest time to the next jump.
We can write
\be \label{eq:recursion1}
T = (R_*+T') \I(R_* < J_*) + (J_*+T_1) \I(J_* < R_*),
\ee
where $T'$ is a random variable with the same distribution as $T$ and independent of $R_*$ and $J_*$, whereas $T_1$ is the end of epidemic time when the two walkers start at a distance one from each other and are both in the infected state. Note that $T_1$ is independent of $R_*$ and $J_*$.
From Equation~\eqref{eq:recursion1}, simple calculations imply that
\be \label{eq:t1}
\L_T(s) = \L_{T_1}(s) \frac{2\lambda}{2\lambda + s},
\ee
where $\L_{T_1}(s) = \E(e^{-sT_1})$ is the Laplace transform of $T_1$.

Similarly, if we define  $R^*\IsDef \max(R_1,R_2)$, we can express $T_1$ as 
\be \label{eq:recursion2}
T_1 = (M+T'') \I(M < R^*) + R^* \I(R^* < M),
\ee
where $T''$ is a random variable with the same distribution as $T$ and independent of $R^*$ and $M$.  Using the independence between the two walkers, the distribution of $R^*$ is simple to obtain and, for every $t\geq 0$, we have $\P(R^* < t) = (1-e^{-\gamma t})^2$, while 
the density is 
\begin{align}
f_{R^*}(t) =
\begin{cases} 2 \gamma e^{-\gamma t}(1-e^{-\gamma t})\;,
  & t \geq 0\;, \\
0\;, & \text{otherwise.}
\end{cases}
\end{align}
Using Equation~\eqref{eq:recursion2} in the Laplace transform of $T_1$, together with the independence between $T''$ and $M$ and $R^*$, we obtain 
\begin{align*}
\L_{T_1}(s) &= \E(e^{-s(M+T'')}; M<R^*) + \E(e^{-sR^*}; R^*<M) \\ &= \L_{T}(s) \E(e^{-sM}; M<R^*) + \E(e^{-sR^*}; R^*<M).
\end{align*}
With the use of Equation~\eqref{eq:t1} to re-write $\L_{T_1}(s)$ on the LHS of the previous equation, we obtain
\be \label{eq:Laplace1}
\L_{T}(s) = \frac{\E(e^{-sR^*}; R^*<M)}{\frac{2 \lambda + s}{2\lambda} - \E(e^{-sM}; M<R^*)}.
\ee 

Let us now deal with the two expectations in Equation~\eqref{eq:Laplace1}. First,
\begin{align*}
\E(e^{-sM}; M<R^*) &= \int_0^\infty e^{-sx} \P(R^*>x) dF_M(x)\\ & = 2\int_0^\infty e^{-sx} e^{-\gamma x} dF_M(x) - \int_0^\infty e^{-sx} e^{-2 \gamma x} dF_M(x)\\ & =  2 \L_M(s+\gamma) - \L_M(s+2 \gamma),
\end{align*}
where $F_M$ is the distribution function of $M$.
Second,
\begin{align*}
&\E(e^{-sR^*}; R^*<M) = \int_0^\infty e^{-sx} \P(M>x)dF_{R^*}(x) = \int_0^\infty e^{-sx} \P(M>x)f_{R^*}(x)dx  \\
&=
2 \gamma \left(\int_0^\infty e^{-(s+\gamma)x} \P(M>x)dx - \int_0^\infty e^{-(s+2\gamma)x} \P(M>x)dx\right)\\ 
& = 
2 \gamma \left(-\int_0^\infty \P(M>x) d\left(\frac{e^{-(s+\gamma)x}}{s+\gamma}\right) + \int_0^\infty \P(M>x) d\left(\frac{e^{-(s+2\gamma)x}}{s+2\gamma}\right)\right)  \\
& =
2\gamma \biggl(\frac{1}{s+\gamma}\left(\left[-\P(M>x) e^{-(s+\gamma)x}\right]_0^\infty - \int_0^\infty e^{-(s+\gamma)x}dF_M(x)\right) \\
& - 
\frac{1}{s+2\gamma}\left(\left[-\P(M>x) e^{-(s+2\gamma)x}\right]_0^\infty - \int_0^\infty e^{-(s+2\gamma)x}dF_M(x)\right) \biggr)  \\
& =
2\gamma \left(\frac{1 - \L_M(s+\gamma)}{s+\gamma} -  \frac{1 - \L_M(s+2\gamma)}{s+2\gamma}\right)
\end{align*}
Combining  everything into Equation~\eqref{eq:Laplace1} we obtain the claim.

\qed

Theorem~\ref{thm:main} provides an expression for the Laplace transform of the end of epidemic time $T$ that is a function of the Laplace transform of the meeting time from distance one $\L_M(s)$, the walking rate $\lambda$ and the recovery rate $\gamma$. This expression determines how the underlying graph structure (encoded in $\L_M(s)$) influences the distribution of $T$. 

\medskip 


The random variable $M$ is related to a discrete time random variable that counts the number of jumps needed for two random walks to meet. More precisely, if $N$ denotes the number of jumps required for the two walkers to meet if they start at distance one, it holds that
\[
M=\sum_{k=1}^N E_k,
\]
where $E_k$ are independent $Exp(2 \lambda)$ random variables.  Note that $N$ only depends on the graph $G$ and does not depend on $\lambda$ nor $\gamma$.  We can then express the Laplace transform of $M$ as
\begin{equation}
\label{eq:MtoN}
\L_M(s) = \sum_{k=1}^\infty \P(N=k) \left(\frac{2 \lambda}{2 \lambda+s}\right)^k = \L_N\left(-\log\left(\frac{2 \lambda}{2 \lambda+s}\right)\right)=\E\left( \left( \frac{2\lambda}{2\lambda + s} \right)^N \right) .
\end{equation}
Combining the above equation with Theorem~\ref{thm:main}, we obtain the following Corollary relating  the Laplace transform of $T$ to the distribution of $N$. 

\begin{corollary}\label{cor:MtoN}
Let $\L_T(s) = \E(e^{-sT})$ denote the Laplace transform of $T$ and let $N$ the number steps before two random walks meet.  Then
\begin{align} \label{eq:Laplace3}
\L_T(s) &= \frac{2\gamma \frac{1 - 
\E\left( \left( \frac{2\lambda}{2\lambda + s+ \gamma} \right)^N \right)
}{s+\gamma} -  \frac{1 - 
\E\left( \left( \frac{2\lambda}{2\lambda + s+2\gamma} \right)^N \right)
}{s+2\gamma}}{\frac{2 \lambda + s}{2\lambda} - 2 
\E\left( \left( \frac{2\lambda}{2\lambda + s+ \gamma} \right)^N \right) 
 + 
 \E\left( \left( \frac{2\lambda}{2\lambda + s+ 2\gamma} \right)^N \right)}
\end{align}
\end{corollary} 

Results of this chapter allow us to completely decouple the effects of the network structure and random walks on it from the effects of the epidemic.

\section{General limit theorems}
\label{sec:limit}


In this section we turn our attention towards understanding how EoE behaves when the size of the underlying graph grows to infinity. To account for the investigation of this asymptotic behavior, let us remark upon two aspects of our  model on  fixed finite graphs:  
\begin{itemize}
\item[$i)$] The distribution of  $T$ may be difficult to explicitly compute;

\item[$ii)$] The  epidemics will eventually die, i.e., $\mathbb{P}(T<+ \infty)=1$ for every value of $\lambda$ and $\gamma$.
\end{itemize}

On the contrary, asymptotically in the size of the graph, we show that there exist scaling regimes for $\lambda$ and $\gamma$  in which the distribution of $T$ can be explicitly computed. Moreover, we also show the existence of regimes in which the epidemic times grow infinitely with the size of the graph.

%

\medskip 

Recall that the random variables $T$, $M$ and $N$ depend on the underlying graph $G$ and, in particular, on its size. With a slight abuse of notation, henceforth we write $T_n$, $M_n$ and $N_n$ to stress the dependence on the graph size $n$, for some  fixed class of graphs (e.g., $K_n$, the complete graph). 

The main objective of this section is to study the limit behavior of $T_n$ when $n$ goes to infinity and to investigate the role of the  underlying graph structure.
To this purpose we  consider walking rate ($\gamma_n$) and recovery rate ($\lambda_n$) which  depend on $n$ such that $\gamma_n\to +\infty$ and $\lambda_n\to +\infty$ as $n$ grows.

This section is organized as follows. In Section \ref{subsec:complete} we consider the complete graph and, by applying Theorem~\ref{thm:main}, we identify limiting behavior of the end of epidemic for all possible scalings of the parameters. Even though Theorem~\ref{thm:main} theoretically allows to deal with any scaling regime and any graph, applying it directly in many cases turns out to be extremely cumbersome. We therefore focus on scalings where the walking dynamics is much faster that the recovery dynamics, so that it is possible that $T_n \to \infty$. We prove two general limit theorems that follow from Theorem~\ref{thm:main} and present their corollaries for specific graphs. Section \ref{subsec:limit_1} contains the first limit theorem followed by its application in Section \ref{subsec:bi} to complete bipartite graphs of certain structures. In Section \ref{subsec:bi} we also identify other bipartite graphs where the first limit theorem is not applicable. In  Section \ref{subsec:limit_2} we prove another limit theorem and discuss its corollaries for bipartite graphs. In Section \ref{subsec:ring} we consider the ring graph, demonstrate that neither limit theorem may be applied to it and study the end of epidemic time in some regimes directly.

\subsection{Complete graph}
\label{subsec:complete}

Let us consider the specific case of a  complete graph on $n$ vertices $G=K_n$.  In this case,  it is not difficult to see that $N_n\sim Geom(n^{-1})$,   $M_n~\sim Exp(\frac{2\lambda}{n})$ and therefore   $\L_{M_n}(s)= 
 \frac{2\lambda_n}{2\lambda_n + n s}$.  Using the latter in Theorem~\ref{thm:main} we obtain that
\begin{align*}
\L_{T_n}(s) &= 
\frac{2\gamma_n \left(\frac{1 - \L_{M_n}(s+\gamma_n)}{s+\gamma_n} -  \frac{1 - \L_{M_n}(s+2\gamma_n)}{s+2\gamma_n}\right)}{\frac{2 \lambda_n + s}{2\lambda_n} - 2 \L_{M_n}(s+\gamma_n) + \L_{M_n}(s+2 \gamma_n)}
\\ & =\frac{2\gamma_n 2\lambda_n \left(\frac{n}{2\lambda_n + n (s+\gamma_n)} -  \frac{n}{2\lambda_n + n (s+2\gamma_n)}\right)}{ 2 \lambda_n + s - 2  \frac{(2\lambda_n)^2}{2\lambda_n + n (s+\gamma_n)} +  \frac{(2\lambda_n)^2}{2\lambda_n + n (s+2\gamma_n)}}
\\
&=
\frac{4\gamma_n^2 \lambda_n  }{ s\left(\frac{2\lambda_n}{n} + s+2\gamma_n\right)\left(\frac{2\lambda_n}{n} +  s+\gamma_n\right)  +(2\lambda_n)^2   \frac{s}{n} + 2\lambda_n (s+2\gamma_n)(s+\gamma_n)}
\\
&= \left( \frac{\lambda_n s}{n^2\gamma_n^2}+ \frac{ s^2}{n\gamma_n^2} + \frac{3 s}{2n\gamma_n} + \frac{s^3}{4 \lambda_n\gamma_n^2} + \frac{3s^2}{4 \lambda_n\gamma_n} + \frac{s}{2 \lambda_n} + \frac{s\lambda_n}{n\gamma_n^2}+ \frac{3s}{2\gamma_n} + \frac{s^2}{2 \gamma_n^2} + 1\right)^{-1}
\end{align*}

Let us consider three different regimes:

\begin{itemize}
\item[i)] $\lambda_n = \omega(n\gamma_n)$ \qquad ({\em i.e., $\frac{\lambda_n}{n \gamma_n}\to +\infty$}) 
\item[ii)] $\lambda_n = o(\gamma_n)$ \qquad ({\em i.e., $\frac{\lambda_n}{\gamma_n}\to 0$})
\item[iii)] $\lambda_n = o(n\gamma_n)$ and $\lambda_n = \omega(\gamma_n)$
\qquad ({\em i.e.,  $\frac{\lambda_n}{\gamma_n}\to + \infty$ but $\frac{\lambda_n}{n \gamma_n}\to 0$}) 
\end{itemize}
\begin{itemize}
\item As far as regime $i)$ is concerned, let us denote by $b_n=\frac{n \gamma_n^2}{\lambda_n}$ and compute $\L_{T_n}(b_n s)$. 

%

\[
\L_{T_n}(b_n s)= \frac{1 } { 1+ s  + o (1) }\; \underset{n\,\uparrow\, \infty}{\longrightarrow}\;\frac{1}{1+s} 
\]

Thus, we have that $ b_n T_n \overset{D}{\longrightarrow}Exp(1)$ as $n$ goes to infinity.
\begin{remark}
Note that within regime $i)$ two different subregimes are possibles:
\begin{itemize}
\item[1)] $b_n \to +\infty$ which implies $T_n \to 0$ in probability,
\item[2)] $b_n \to 0$ which implies  $T_n \to + \infty$ in probability.
\end{itemize}
The regime $2)$ is of particular interest as the epidemic time grows infinitely, despite being finite for every $n$.
\end{remark}
Note that $b_n \to 0$ implies that $\lambda_n = \omega(n\gamma_n^2)$ and, thus the walkers are moving at a rate that scales linearly with $n$ and quadratically with the recovery rate. Interestingly, this rate is large enough for the epidemic time to grow infinitely. 

\item As far as regime $ii)$ is concerned, let us  compute $\L_{T_n}(2\lambda_n  s)$. 
%
\[
\L_{T_n}(2\lambda_n s)=\frac{1 } { 1+ s  + o (1) }\; \underset{n\,\uparrow\, \infty}{\longrightarrow}\;\frac{1}{1+s} 
\]
Thus, $ 2\lambda_n  T_n \overset{D}{\longrightarrow}Exp(1)$ as $n$ goes to infinity.

\item As far as regime $iii)$ is concerned, let us  compute $\L_{T_n}(\gamma_n  s)$
\[
\L_{T_n}(\gamma_n s)= \frac{1}{ ( s  + 1)(s/2 +1) +  o(1) } \underset{n\uparrow \infty}{\longrightarrow}\frac{2}{(s+1)(s +2)}\;.
\]
Thus, $ \gamma_n  T_n \overset{D}{\longrightarrow} X + Y$ as $n$ goes to infinity, where $X\sim Exp(1)$ $Y\sim Exp(2)$ and $X$ and $Y$ are independent.
\end{itemize}
 
\begin{remark}
In both regimes $ii)$ and $iii)$  we have that $T_n\to 0$ in probability.
\end{remark}

\medskip

\subsection{General limit theorem I}
\label{subsec:limit_1}

For a complete graph, we have shown how Theorem~\ref{thm:main} can be applied to identify  the behaviour of $T_n$ in different regimes. This same approach can be applied to other classes of graphs. However, instead of carrying out a similar analysis for other graphs, we provide below (and also in Section \ref{subsec:limit_2}) auxiliary results which allow us to characterize the asymptotic behaviour of $T_n$ in terms of the asymptotic behaviour of $N_n$. These results can then be directly applied to other classes of graphs.

\begin{theorem}\label{thm:limit_general}
Assume  there exists a sequence $\{a_n\}_{n \in \mathbb{N}}$ converging to zero such that  $a_n N_n \overset{D}{\longrightarrow}X$ as $n$ tends to infinity, with $X$ a random variable with the first two moments finite. Denote $c_1 = \E(X)$ and $c_2 = \E(X^2)$. Then, for $\gamma_n = o(\lambda_n a_n)$ 
\[
\frac{c_2}{2c_1}\frac{ \gamma_n^2 }{\lambda_n a_n}T_n \overset{D}{\longrightarrow} \text{Exp}(1)\;.
\] 

\end{theorem}

%

%

The above theorem easily captures the regime i) for the complete graph. Indeed,  for the complete graph,  $N_n \sim Geom (\frac{1}{n})$ and  $\L_{N_n}(s)= \frac{e^{-s}}{n - (n-1)e^{-s}}$. Thus, 
\[
\L_{N_n}\left(n^{-1}s\right)= \frac{e^{-\frac{s}{n}}}{n - (n-1)e^{-\frac{s}{n}}}= \frac{1 - \frac{s}{n} + o(\frac{1}{n})}{s + o(1) + 1 - \frac{s}{n}} \underset{n\uparrow \infty}{\longrightarrow}=\frac{1}{1+s}\;.
\]  
Therefore,  $n^{-1}N_n\overset{D}{\longrightarrow}Exp(1)$ as $n$ tends to infinity. In Theorem~\ref{thm:limit_general} this corresponds to  the situation where $X\sim Exp(1)$ and  thus $c_1=1$, $c_2=2$, and $a_n = n^{-1}$.   In the regime $\frac{\lambda_n}{n \gamma_n}\to +\infty$, the above theorem guarantees that  $ \frac{n \gamma_n^2}{\lambda_n} T_n \overset{D}{\longrightarrow}Exp(1)$ as $n$ grows, as previously shown through several calculations.

\begin{remark}
Let us give an intuitive explanation of why the regime $\gamma_n = o(\lambda_n a_n)$ is interesting. Due to the assumption $a_n N_n \overset{D}{\longrightarrow}X$, it takes on average order $1/a_n$ steps for the two random walks to meet as soon as they separate. Each step takes an average time of order $1/\lambda_n$, so it will take on average time of order $1/(\lambda_n a_n)$ for the two walkers to meet. On the other hand, it will take an average time of order $1/\gamma_n$ for them both to recover. Therefore, assuming $\gamma_n = o(\lambda_n a_n)$, we make sure that the time for the two walkers to recover is much larger than the time it takes them to meet again, thus giving the epidemic a chance to survive.
\end{remark}

{\bf Proof of Theorem~\ref{thm:limit_general}.}

To simplify the notation let us denote $b_n=\dfrac{c_2 \gamma_n^2}{2c_1\lambda_n a_n}$.  To prove the claim it is enough to show that
$\L_{T_n}(b_n s)\underset{n \uparrow \infty}{\longrightarrow}\frac{1}{1+s}$, for every $s>0$. 
The assumption $a_n N_n \overset{D}{\longrightarrow}X$ assures that 
$\L_{N_n}(a_n s) \underset{n \uparrow \infty}{\longrightarrow} \L_X(s)$.
From the latter   we obtain a limit theorem for $M_n$, i.e., 
\begin{align*}
\L_{M_n}(2\lambda_n a_n s) & = \L_{N_n}\left(-\log\left(\frac{1}{1+a_n s}\right)\right)
\\ &= \L_{N_n}\left(\log \left(1+a_n s\right) \right)  = \L_{N_n}\big(a_n s + o(a_n s)\big) \underset{n\uparrow \infty}{\longrightarrow} \L_X(s)\;,
\end{align*}
where we used  that $a_n$ goes to zero and the continuity of $\L_X$. 

\medskip 

For $s$ sufficiently close to zero we can write  $\L_X(s)=1 - c_1 s + \frac{c_2}{2} s^2  + o(s^2)$. Thus, for  any sequence 
$x_n \to 0$ it holds that 
$\L_X(x_ 
n s) = 1-c_1x_n s + \frac{c_2}{2}(x_n s)^2 + o ((x_n s)^2)$.

We use Corollary~\ref{cor:MtoN} to compute $\L_{T_n}(b_n s)$ from $\L_{N_n}$. We begin computing  $\L_{M_n}\left(b_n s+\gamma_n\right)$.
\begin{align*}
\L_{M_n}\left(b_n s +\gamma_n\right) &= 
\L_{M_n}\left(2\lambda_n a_n\frac{b_n s+\gamma_n}{2\lambda_n a_n }\right) = 
 \L_{N_n}\left(-\log\left(\frac{1}{1+a_n \frac{b_n s+\gamma_n}{2\lambda_n a_n}}\right)\right)
 \\&
 =\L_{N_n}\left(a_n \frac{ b_n s+\gamma_n}{2\lambda_n a_n} + o\left(a_n \frac{b_n s+\gamma_n}{2\lambda_n a_n }\right)\right)
\end{align*}
Using that $\frac{1}{\lambda_n a_n}\to 0$, we have that 
\begin{align*}
\L_{N_n}\left(a_n \frac{ b_n s+\gamma_n}{2\lambda_n a_n} + o\left( a_n \frac{ b_n s+\gamma_n}{2\lambda_n a_n} \right)\right) \\ = 1 - c_1 \left(\frac{b_n s + \gamma_n }{2\lambda_n a_n} \right) + \frac{c_2}{2} \left(\frac{b_n s + \gamma_n }{2\lambda_n a_n} \right)^2  + o\left(\left(\frac{\gamma_n}{\lambda_n a_n} \right)^2 \right).
\end{align*}
So we obtain
$\L_{M_n}\left(b_n s+\gamma_n\right)= 1 - c_1 \left(\frac{b_n s + \gamma_n}{2\lambda_n a_n } \right) + \frac{c_2}{2} \left(\frac{b_n s + \gamma_n}{2\lambda_n a_n } \right)^2  + o\left(\left(\frac{\gamma_n}{\lambda_n a_n} \right)^2\right)
$. 
Similarly, $\L_{M_n}\left(b_n s+2\gamma_n\right)= 1 - c_1 \left(\frac{b_n s +2\gamma_n}{2\lambda_n a_n } \right) + \frac{c_2}{2} \left(\frac{b_n s+2\gamma_n}{2\lambda_n a_n } \right)^2  + o\left(\left(\frac{\gamma_n }{\lambda_n a_n} \right)^2\right)
$. Therefore, the numerator in Equation~\eqref{eq:Laplace2} can be written as 
\begin{align*}
&2\gamma_n \left(\frac{1 - \L_{M_n}(b_n s+\gamma_n)}{b_n s+\gamma_n} -  \frac{1 - \L_{M_n}(b_n s+2\gamma_n)}{b_n s+2\gamma_n}\right)\\
& = \frac{2\gamma_n}{(b_n s+\gamma_n)} \left( c_1 \left(\frac{b_n s + \gamma_n}{2\lambda_n a_n } \right) - \frac{c_2}{2} \left(\frac{b_n s + \gamma_n}{2\lambda_n a_n } \right)^2  + o\left( \left(\frac{\gamma_n}{\lambda_n a_n } \right)^2\right) \right)
\\
& - 
\frac{2\gamma_n}{(b_n s+2\gamma_n)} \left( c_1 \left(\frac{b_n s + 2\gamma_n}{2\lambda_n a_n } \right) - \frac{c_2}{2} \left(\frac{b_n s + 2\gamma_n}{2\lambda_n a_n } \right)^2  + o\left( \left(\frac{\gamma_n}{\lambda_n a_n } \right)^2\right) \right)
\\
&=
 c_1 \left(\frac{\gamma_n}{\lambda_n a_n } \right) - c_2 \left(\frac{\gamma_n}{2\lambda_n a_n } \right)^2  + o\left( \left(\frac{\gamma_n}{\lambda_n a_n } \right)^2\right) 
\\
& - c_1 \left(\frac{\gamma_n}{\lambda_n a_n } \right) + 2 c_2 \left(\frac{\gamma_n}{2\lambda_n a_n } \right)^2  + o\left( \left(\frac{\gamma_n}{\lambda_n a_n } \right)^2\right) 
\\
&= c_2 \left(\frac{\gamma_n}{2\lambda_n a_n } \right)^2  + o\left( \left(\frac{\gamma_n}{\lambda_n a_n } \right)^2\right)\;.
\end{align*}
Let us now look at the denominator:
\begin{align*}
&\frac{2 \lambda_n + b_n s}{2\lambda_n} - 2 \L_M(b_n s+\gamma_n) + \L_M(b_n s+2 \gamma_n) = \\
&= 1 + \frac{b_n s}{2 \lambda_n}  -2\left( 1 - c_1 \left(\frac{b_n s + \gamma_n}{2\lambda_n a_n } \right) + \frac{c_2}{2} \left(\frac{b_n s + \gamma_n}{2\lambda_n a_n } \right)^2  + o\left(\left(\frac{\gamma_n}{\lambda_n a_n} \right)^2\right) \right) 
\\
&+  1 - c_1 \left(\frac{b_n s+ 2\gamma_n}{2\lambda_n a_n } \right) + \frac{c_2}{2} \left(\frac{b_n s+ 2\gamma_n}{2\lambda_n a_n } \right)^2  + o\left(\left(\frac{\gamma_n}{\lambda_n a_n} \right)^2\right)
\\
&= \frac{b_n s}{2 \lambda_n} + c_1  \left(\frac{b_n s}{2\lambda_n a_n }  \right) + c_2 \left(\frac{\gamma_n}{2\lambda_n a_n} \right)^2  +  o\left(\left(\frac{\gamma_n}{\lambda_n a_n} \right)^2\right)
\end{align*}

Overall, recalling that $b_n= \frac{c_2\gamma_n^2}{2c_1\lambda_n a_n}$, we have

\begin{align*}
\L_{T_n}(b_n s) & = \frac{c_2 \left(\frac{\gamma_n}{2\lambda_n a_n } \right)^2  + o\left( \left(\frac{\gamma_n}{\lambda_n a_n } \right)^2\right) }{\frac{b_n s}{2 \lambda_n} + c_1  \left(\frac{b_n s}{2\lambda_n a_n }  \right) + c_2 \left(\frac{\gamma_n}{2\lambda_n a_n} \right)^2  +  o\left(\left(\frac{\gamma_n}{\lambda_n a_n} \right)^2\right)} 
\\ &= \frac{1+ o(1)}{s \frac{a_n}{c_1} + s + 1 + o(1)} \underset{n \uparrow \infty}{\longrightarrow}  \frac{1}{1+ s}
\end{align*}

\qed

We have already seen how the above theorem can be applied to the complete graph in the regime $\gamma_n = o (\lambda_n n)$. 
Another example in which Theorem~\ref{thm:limit_general} can be applied is the complete bipartite graph $K_{m,n-m}$ presented below.

\subsection{Complete bipartite graph}\label{sec:complete-bipartite}
\label{subsec:bi}
For the complete bipartite graph $G=K_{m,n-m}$ on $n$ vertices with one partition having $m$ elements and the other $n-m$ elements, we shall consider  several different scenarios:
 i)  $m=\alpha n$, with  $\alpha \in (0,1)$,  ii) $m=n^\beta$, with $\beta \in (0,1)$ and iii) $m=\log n$.
%
%
It is not difficult to see that, in this case the random variable $N_n$ satisfies
\begin{eqnarray*}
N_n \overset{d}{=}
\begin{cases}
1, \quad \text{w.p.} \quad \frac{1}{2m}, \\
1, \quad \text{w.p.} \quad  \frac{1}{2(n-m)}, \\
2+N_n, \quad \text{w.p.} \quad 1- \frac{1}{2m} -\frac{1}{2(n-m)}.
\end{cases}
\end{eqnarray*}
Therefore the Laplace transform of $N_n$  satisfies he following recursion 
$$
\L_{N_n}(s) = \frac{1}{2} \left(\frac{1}{m} + \frac{1}{n-m} \right) e^{-s} + \frac{1}{2} \left( 2 - \frac{1}{m} - \frac{1}{n-m}  \right)e^{-2s}\L_{N_n}(s)\;,
$$
and hence
$$
\L_{N_n}(s) = \frac{\frac{1}{2} \left(\frac{1}{m} + \frac{1}{n-m}\right) e^{-s}}{ 1- \frac{1}{2}\left( 2 - \frac{1}{m} - \frac{1}{n-m}  \right)e^{-2s}}= \frac{\frac{1}{2} \left(\frac{n}{m(n-m)} \right) e^{-s}}{ 1- e^{-2s} +  \frac{1}{2}\left( \frac{n}{m(n-m)}  \right)e^{-2s}} \;.
$$ 

\begin{itemize}
%
\item If $m=\alpha n$ we have
\[
\L_{N_n}(s)= 
 \frac{\frac{1}{2\alpha(1-\alpha)}e^{-s}}{n(1-e^{-2s}) + \frac{1}{2\alpha(1-\alpha)}e^{-s}}\;,
\]
and therefore,  
\begin{align*}
\L_{N_n}(n^{-1}s)&=  \frac{\frac{1}{2\alpha(1-\alpha)}e^{-\frac{s}{n}}}{n(1-e^{-\frac{2s}{n}}) + \frac{1}{2\alpha(1-\alpha)}e^{-\frac{s}{n}}}
\\ &= \frac{\frac{1}{2\alpha(1-\alpha)}(1 - \frac{s}{n} + o(\frac{1}{n}))}{n(1-(1 - 2\frac{s}{n} + o(\frac{1}{n}))) + \frac{1}{2\alpha(1-\alpha)}(1 - \frac{s}{n} + o(\frac{1}{n}))}
\\
&=
\frac{\frac{1}{2\alpha(1-\alpha)} (1 + o(1))}{ 2 s +  \frac{1}{2\alpha(1-\alpha)} +  o(1)} \underset{n\uparrow \infty}{\longrightarrow} 
\frac{(4\alpha(1-\alpha))^{-1}}{(4\alpha(1-\alpha))^{-1} +  s}\;,
\end{align*}
that is $n^{-1}N_n \overset{D}{\longrightarrow}Exp ((4\alpha(1-\alpha))^{-1})$. 
Thus, Theorem~\ref{thm:limit_general} guarantees that $\frac{c_2}{2c_1}\frac{ \gamma_n^2 }{\lambda_n a_n}T_n \overset{D}{\longrightarrow} \text{Exp}(1)$, with $c_1=4\alpha(1-\alpha)$ and $c_2=32\alpha^2(1-\alpha)^2$.

\item If $m=m(n)\uparrow + \infty$ and $m=o(n)$

\[
\L_{N_n}(s)= 
\frac{\frac{1}{2} \left(\frac{n}{m(n-m)} \right) e^{-s}}{ 1- e^{-2s} +  \frac{1}{2}\left( \frac{n}{m(n-m)}  \right)e^{-2s}}=\frac{ \left(\frac{1}{1-o(1)} \right) e^{-s}}{ 2m(1- e^{-2s}) +  \left( \frac{1}{1-o(1)}  \right)e^{-2s}}
\]
and therefore, 
\begin{align*}
\L_{N_n}(m^{-1}s)&=  
\frac{\left(\frac{1}{1-o(1)} \right)\left(1 - \frac{s}{m} + o(\frac{1}{m})\right)}{2m ( 2\frac{s}{m} + o(\frac{1}{m}))) + \left(\frac{1}{1-o(1)} \right)\left(1 - \frac{s}{m} + o(\frac{1}{m})\right)}
\underset{n\uparrow \infty}{\longrightarrow} 
\frac{4^{-1}}{4^{-1} +  s}
\end{align*}

that is, $m^{-1}N_n \overset{D}{\longrightarrow}Exp (4^{-1})$  and 
Theorem~\ref{thm:limit_general} implies  that 
$$
\frac{c_2}{2c_1}\frac{ \gamma_n^2 }{\lambda_n m^{-1}}T_n \overset{D}{\longrightarrow} \text{Exp}(1),
$$ 
with $c_1=4$ and $c_2=32$.

Two notable examples of this scenarios are:
\begin{itemize}
\item[i)]  {\em power law growth}: $m(n)=n^\beta$, with $\beta \in (0,1)$
\item[ii)] {\em polylogarithmic growth}:  $m(n)=\log^\beta n$, with $\beta>0$
\end{itemize}

\end{itemize}


\medskip

There are situations  in which Theorem~\ref{thm:limit_general} cannot be applied. Specifically, it is not always the case that, given $a_n$ going to zero,  $a_n N_n$  converges in distribution to a random variable.  Consider, for  example, the complete bipartite graph $K_{m,n-m}$ with $m$ constant. 
 In this case, the Laplace transform is given by 
\[
\L_{N_n}(s) = \frac{\frac{1}{2} \left(\frac{1}{m} + \frac{1}{n-m}\right) e^{-s}}{ 1- \frac{1}{2}\left( 2 - \frac{1}{m} - \frac{1}{n-m}  \right)e^{-2s}}= \frac{\frac{1}{2} \left(\frac{n}{m(n-m)} \right) e^{-s}}{ 1- e^{-2s} +  \frac{1}{2}\left( \frac{n}{m(n-m)}  \right)e^{-2s}} \;.
\] 
and, as $m$ is a constant, we have that
$
\L_{N_n}(s) \underset{n\uparrow \infty}{\longrightarrow}
\frac{\frac{1}{2m}  e^{-s}}{ 1- \frac{1}{2}\left( 2 - \frac{1}{m}\right)e^{-2s}}
$.
%
%
%
%
%
As it turns out,  in this case $N_n$  converges in distribution to a random variable with first and second moments finite, but  $a_n N_n \to 0$ in probability for any  sequence $a_n \to 0$.
In the next section we provide a general result to characterize the limit behaviour of $T_n$  in these type of situations. 

\subsection{General limit theorem II}
\label{subsec:limit_2}

We present a limit theorem in a different regime to that of Theorem \ref{thm:limit_general}. This will allow to consider the bipartite graphs identified in the previous section for which Theorem \ref{thm:limit_general} is not applicable.

\begin{theorem}\label{thm:limit_general_2}
Assume  $N_n \overset{D}{\longrightarrow}X$ as $n$ tends to infinity, with $X$ a random variable with the first two moments finite. Denote $c_1 = \E(X)$ and $c_2 = \E(X^2)$. Then, for $\gamma_n = o(\lambda_n)$ 
\[
\frac{(c_1 + c_2)}{2(1+c_1)}\frac{ \gamma_n^2 }{\lambda_n}T_n \overset{D}{\longrightarrow} \text{Exp}(1)\;.
\] 
\end{theorem}

Note that, for  the complete bipartite  graph with a finite fixed $m$, we have that $\L_{N_n}(s)$ converges to 
$\frac{1/(2m) e^{-s}}{1-1/2(2-1/m)e^{-2s}}$. Thus, $N_n$ converges in distribution to a random variable  $X$ with $c_1=\E(X)=-\frac{d}{ds}\L_X(s)\vert_{s=0}=4m-1 $  and $c_2=\E(X^2)=\frac{d^2}{ds^2}\L_X(s)\vert_{s=0}= 16m(2m-1)+1$. Thus the above theorem guarantees that, in the regime $\gamma_n = o(\lambda_n)$, the random variable  $\frac{8m -3}{2}\frac{ \gamma_n^2 }{\lambda_n}T_n $ converges in distribution to an exponential random variable.
In the particular case  when $m=1$, which corresponds  to the star graph,  we have that $\frac{5}{2}\frac{ \gamma_n^2 }{\lambda_n}T_n $ converges in distribution to an exponential random variable.

{\bf Proof of Theorem~\ref{thm:limit_general_2}.}

The proof follows the same approach as the proof of Theorem~\ref{thm:limit_general}, and is provided for completeness. 
To simplify the notation let us define $b_n=\frac{(c_1 + c_2)}{2(1+c_1)}\frac{ \gamma_n^2 }{\lambda_n}$.  To show the claim it is enough to show that
$\L_{T_n}(b_n s)\underset{n \uparrow \infty}{\longrightarrow}\frac{1}{1+s}$, for every $s>0$. 
The assumption $N_n \overset{D}{\longrightarrow}X$ assures that 
$\L_{N_n}(s) \underset{n \uparrow \infty}{\longrightarrow} \L_X(s)$, where $\L_X(s)=1 - c_1 s + \frac{c_2}{2} s^2  + o(s^2)$.
From the latter, using Equation~\ref{eq:MtoN}, we obtain a limit theorem for $M_n$, i.e., 
\[
\L_{M_n}(2\lambda_n  s) = \L_{N_n}\left(-\log\left(\frac{1}{1+s}\right)\right) = \L_{N_n}(\log(1+s)) \underset{n\uparrow\infty}{\longrightarrow} \L_X(\log(1+s))\;.
\]

Using the Taylor expansions for $\L_X$ and for $\log$, we can write
\begin{align*}
\L_X(\log(1+s)) & = 1 - c_1\log(1+s) + \frac{c_2}{2} (\log(1+s))^2 + o(s^2) \\ &= 1 - c_1 \big(s - s^2/2+o(s^2)\big) + \frac{c_2}{2}\big(s-s^2/2 + o(s^2)\big)^2 + o(s^2) \\ &= 1 - c_1s +\frac{c_1 + c_2}{2} s^2 + o(s^2)\;.
\end{align*}

Thus,
\begin{align*}
\L_{M_n}\left(b_n s+\gamma_n\right) &= \L_{M_n}\left(2\lambda_n \frac{b_n s+\gamma_n}{2\lambda_n }\right) = 
 \L_{N_n}\left(\log\left(1+ \frac{b_n s+\gamma_n}{2\lambda_n }\right)\right)
 \\&
 = 1 - c_1 \left(\frac{b_n s + \gamma_n }{2\lambda_n } \right) + \frac{c_1+c_2}{2} \left(\frac{b_n s + \gamma_n }{2\lambda_n } \right)^2  + o\left(\left(\frac{\gamma_n}{\lambda_n } \right)^2 \right)
\end{align*}
Similarly, $\L_{M_n}\left(b_n s+2\gamma_n\right)= 1 - c_1 \left(\frac{b_n s +2\gamma_n}{2\lambda_n } \right) + \frac{c_1+c_2}{2} \left(\frac{b_n s+2\gamma_n}{2\lambda_n } \right)^2  + o\left(\left(\frac{\gamma_n }{\lambda_n } \right)^2\right)
$. Therefore, the numerator in Equation~\eqref{eq:Laplace2} can be written as 
\begin{align*}
&2\gamma_n \left(\frac{1 - \L_{M_n}(b_n s+\gamma_n)}{b_n s+\gamma_n} -  \frac{1 - \L_{M_n}(b_n s+2\gamma_n)}{b_n s+2\gamma_n}\right)\\
& = \frac{2\gamma_n}{(b_n s+\gamma_n)} \left( c_1 \left(\frac{b_n s + \gamma_n}{2\lambda_n } \right) - \frac{c_1+c_2}{2} \left(\frac{b_n s + \gamma_n}{2\lambda_n } \right)^2  + o\left( \left(\frac{\gamma_n}{\lambda_n  } \right)^2\right) \right)
\\
& - 
\frac{2\gamma_n}{(b_n s+2\gamma_n)} \left( c_1 \left(\frac{b_n s + 2\gamma_n}{2\lambda_n } \right) - \frac{c_1+c_2}{2} \left(\frac{b_n s + 2\gamma_n}{2\lambda_n } \right)^2  + o\left( \left(\frac{\gamma_n}{\lambda_n  } \right)^2\right) \right)
\\
&=
 c_1 \left(\frac{\gamma_n}{\lambda_n  } \right) - (c_1+c_2) \left(\frac{\gamma_n}{2\lambda_n } \right)^2  + o\left( \left(\frac{\gamma_n}{\lambda_n } \right)^2\right) 
\\
& - c_1 \left(\frac{\gamma_n}{\lambda_n  } \right) + 2 (c_1+c_2) \left(\frac{\gamma_n}{2\lambda_n } \right)^2  + o\left( \left(\frac{\gamma_n}{\lambda_n } \right)^2\right) 
\\
&= (c_1+c_2) \left(\frac{\gamma_n}{2\lambda_n } \right)^2  + o\left( \left(\frac{\gamma_n}{\lambda_n } \right)^2\right)\;.
\end{align*}
Let us now look at the denominator:
\begin{align*}
&\frac{2 \lambda_n + b_n s}{2\lambda_n} - 2 \L_M(b_n s+\gamma_n) + \L_M(b_n s+2 \gamma_n) = \\
&= 1 + \frac{b_n s}{2 \lambda_n}  -2\left( 1 - c_1 \left(\frac{b_n s + \gamma_n}{2\lambda_n } \right) + \frac{c_1+c_2}{2} \left(\frac{b_n s + \gamma_n}{2\lambda_n } \right)^2  + o\left(\left(\frac{\gamma_n}{\lambda_n } \right)^2\right) \right) 
\\
&+  1 - c_1 \left(\frac{b_n s+ 2\gamma_n}{2\lambda_n } \right) + \frac{c_1+c_2}{2} \left(\frac{b_n s+ 2\gamma_n}{2\lambda_n  } \right)^2  + o\left(\left(\frac{\gamma_n}{\lambda_n } \right)^2\right)
\\
&= \frac{b_n s}{2 \lambda_n} + c_1  \left(\frac{b_n s}{2\lambda_n }  \right) + (c_1+c_2) \left(\frac{\gamma_n}{2\lambda_n} \right)^2  +  o\left(\left(\frac{\gamma_n}{\lambda_n } \right)^2\right)
\end{align*}
Overall, recalling that $b_n= \frac{(c_1+c_2)\gamma_n^2}{2(1+c_1)\lambda_n }$, we have
\begin{align*}
\L_{T_n}(b_n s) = \frac{(c_1+c_2) \left(\frac{\gamma_n}{2\lambda_n  } \right)^2  + o\left( \left(\frac{\gamma_n}{\lambda_n } \right)^2\right) }{\frac{b_n s}{2 \lambda_n}(1+  c_1)    + (c_1+c_2) \left(\frac{\gamma_n}{2\lambda_n } \right)^2  +  o\left(\left(\frac{\gamma_n}{\lambda_n } \right)^2\right)}&= \frac{1+ o(1)}{s + 1 + o(1)} \\
 &\underset{n \uparrow \infty}{\longrightarrow}  \frac{1}{1+ s}
\end{align*}

\qed

\subsection{Ring}
\label{subsec:ring}

In this section, we study the behaviour of $T_n$ on the ring $C_n$ on $n$ vertices. In the sequel, we shall assume that $n$ is even; this is not crucial and the case $n$ odd, albeit slightly different, can be similarly handled. 

Following the path already used in the previous sections, we need to first understand the behaviour of the number of steps two walkers need to meet up starting at distance one. In light of this, for $i=1, \ldots, \frac{n}{2}$, let $N_{n,i}$ be the number of steps two random walks need to meet starting at distance $i$ on the ring on $n$ vertices; we set $N_{n,0}=0$. 
 The following recursion holds
\begin{align}\label{eq:recur}
N_{n,i}= \begin{cases}
1 + N_{n,i-1} & \text{w.p. $\frac{1}{2}$}
\\
1 + N_{n,i+1} & \text{w.p. $\frac{1}{2}$}
\end{cases}
\;\;, \text{ for $1\leq i< \frac{n}{2}$}
\end{align}
and
$$
N_{n,\frac{n}{2}}=
1 + N_{n,\frac{n}{2}-1} \qquad  \text{ w.p. $1$.}
$$
To simplify the notation let us set $\L_i(s)=\L_{N_{n,i}}(s)$ to denote the Laplace transform of $N_{n,i}$, with $\L_0(s)=1$, and  $\alpha=e^{-s}/2$. Using Equation~\ref{eq:recur}  we obtain the following recursion for $\L_i$
\begin{align*}
\L_i (s)&= \alpha \L_{i-1}(s) + \alpha \L_{i+1}(s) \;  \qquad \text{ for $1\leq i \leq \frac{n}{2}-1$}
\\  
\L_{\frac{n}{2}}(s) &= 2 \alpha \L_{\frac{n}{2}-1}(s)
\end{align*}
Recall that we are interested in $\L_1(s)$, i.e., the Laplace transform of the number of steps two walkers need to meet starting at distance one.
Solving the  latter recursion  (see Appendix~\ref{app:recursion1} for details) we obtain 
\begin{align}\label{eq:ring}
\L_1(s)=\frac{\alpha}{x_1}\; \frac{1  + \left( \frac{x_2}{x_1} \right)^{\frac{n}{2}-1}   }{ 1  + \left( \frac{x_2}{x_1} \right)^{\frac{n}{2}} }\;,
&& \L_1(s) \underset{n\uparrow\infty}{\longrightarrow} \frac{\alpha}{x_1}=\frac{e^{-s}}{1+\sqrt{1-e^{-2s}}}
\end{align}
where $x_1= \frac{1 + \sqrt{1-4\alpha^2}}{2}$ and $x_2= \frac{1 - \sqrt{1-4\alpha^2}}{2}$.

\medskip 

From Equation~\ref{eq:ring} two observations can be made:
\begin{itemize}
     \item $a_n N_{n,1} \overset{\mathbb{P}}{\longrightarrow} 0$ for all $a_n \to 0$ which  implies that Theorem~\ref{thm:limit_general} cannot be applied 
    \item $N_n\overset{D}{\longrightarrow}X$ with  $\mathbb{E}(X)=+\infty$ which implies that   Theorem~\ref{thm:limit_general_2} cannot be applied 
\end{itemize}
\begin{remark}
Note that the second observation makes sense 
as $N_{n,1}$ should converge to the time for a simple symmetric random walk on integers to reach zero if it starts at $1$, which is well-known to have an infinite expectation.
\end{remark} 

In order to study the asymptotic of $T_n$ on the ring, we therefore resort to Theorem~\ref{thm:main} and obtain 
\[
\L_{M_n}(s) = \L_{N_{n,1}}\left(\log\left(\frac{2 \lambda_n}{2 \lambda_n+s}\right)^{-1}\right) 
=
\frac{2\lambda_n}{2\lambda_n+s+\sqrt{s(4\lambda_n+s)}} + o(1)\;.
\]
Let us restrict to the regime $\gamma_n = o(\lambda_n)$ and compute $\L_{T_n}(\gamma_n s)$. First we observe that
\begin{align*}
   1-\L_{M_n}(\gamma_n (1+s)) &= \frac{\sqrt{\gamma_n(1+s)}}{\sqrt{\lambda_n}} (1+o(1)) 
   \\
 1-\L_{M_n}(\gamma_n (2+s))  &= \frac{\sqrt{\gamma_n(2+s)}}{\sqrt{\lambda_n}} (1+o(1))
 \end{align*}
Thus, 
 \begin{align*}
 \L_{T_n}&(\gamma_n s)= \frac{2\gamma_n \left(\frac{1 - \L_{M_n}(\gamma_n(s+1))}{\gamma_n(s+1)} -  \frac{1 - \L_{M_n}((s+2)\gamma_n)}{(s+2)\gamma_n}\right) }{ 1 + \frac{\gamma_n s}{2 \lambda_n} - 2 + 2(1-\L_{M_n}(\gamma_n s)) + 1 - (1 - \L_{M_n}(\gamma_n s))}
 \\
 &=
 \frac{\frac{2\sqrt{\gamma_n}(\sqrt{2+s} - \sqrt{1+s})}{\sqrt{\lambda_n} \sqrt{1+s} \sqrt{2+s}}  (1+o(1))}{ \frac{\sqrt{\gamma_n}(2\sqrt{1+s} - \sqrt{2+s})}{\sqrt{\lambda_n}}(1+o(1))}\; \longrightarrow \; \frac{2(\sqrt{2+s} - \sqrt{1+s})}{\sqrt{1+s} \sqrt{2+s}(2\sqrt{1+s} - \sqrt{2+s})}\;.
\end{align*}

The latter tells us that $\gamma_n T_n \overset{D}{\longrightarrow} X$ where
$X$ is a random variables such that 
$\L_X(s)= \frac{2(\sqrt{2+s} - \sqrt{1+s})}{\sqrt{1+s} \sqrt{2+s}(2\sqrt{1+s} - \sqrt{2+s})}$.


\section{Final Remarks}
\label{sec:conclusion}

Epidemics on networks driven by mobile agents serve as a fundamental model for different contagious processes, finding applications in various domains. In the SIS epidemic model agents alternate between being susceptible and infected, becoming infected when meeting in network nodes, and a fundamental statistic is the duration of the epidemic (since all agents will eventually become susceptible). This model is challenging to analyze due the dependence between the epidemic process and agent mobility. When agent mobility is agnostic to the epidemic process (e.g., agents perform independent random walks, the scenario tackled in this work), theoretical analysis is more manageable. 

Indeed, by considering edge-transitive graphs and two agents, this work establishes a strong result that separates the epidemic process from the meeting process. In particular, Theorem~\ref{thm:main} determines the Laplace transform of the epidemic end time (EoE) as a function of the Laplace transform of the meeting times. Note that the latter depends only on the network structure. The second contribution is the characterization of the EoE for graph sequences of increasing size (Theorems~\ref{thm:limit_general} and~\ref{thm:limit_general_2}). While for every finite graph, the EoE is finite, under a proper scaling of the model parameters the EoE can be arbitrarily long (and even converge to infinity, as the graph size grows). Interestingly, the proper scaling for such phenomenon strongly depend on the graph structure. This finding highlights a possible phase transition between very short and very long (expected) EoE, a phenomenon that has been rigorously observed in a related model~\cite{dickman2010activated,sidoravicius2017absorbing}.

While this work focused on two agents, a natural next step is the characterization of EoE as a function of the number of agents. Indeed, recent works on a related model have shown that the density of the number of agents (in infinite lattices of fixed dimension), plays a fundamental role on EoE~\cite{dickman2010activated,sidoravicius2017absorbing}. While the approach taken in this work does not trivially extend to three agents, a mean-field approach could be derived for finite graphs with a sufficient number of agents. In fact, we conjecture that with a large enough number of agents, the EOE will be similar to the EoE in the classic network epidemic model, where network nodes have epidemic states. This result would establish an important relationship between apparently different models, contributing further to our understanding of network epidemics.

\appendix
\section{Appendix}\label{app:recursion1}
Hereby we solve the recursion 
\begin{align*}
\L_i (s)&= \alpha \L_{i-1}(s) + \alpha \L_{i+1}(s) \;  \qquad \text{ for $1\leq i \leq \frac{n}{2}-1$}
\\  
\L_{\frac{n}{2}}(s) &= 2 \alpha \L_{\frac{n}{2}-1}(s)
\end{align*}
Let latter expression can be rewritten in the following 
form
\begin{align*}
\L_i(s)= C_i \L_{i-1}(s)\;, \qquad  \text{ for $1\leq i \leq \frac{n}{2}$}
\end{align*}
where, for instance,  $C_{\frac{n}{2}}=2\alpha$ and $C_{\frac{n}{2}-1}= \frac{\alpha}{1-2\alpha^2}$.
Note that, given that $\L_0(s)=1$, we have  $\L_1(s)=C_1$. Thus, the problem of finding the Laplace transform of $N_{n,1}$ reduces  to compute $C_1$.
The coefficients $C_i$ satisfy the following recursion 
\begin{align}\label{eq:C2}
C_i = \frac{\alpha}{(1- \alpha C_{i+1})}\;, \quad \text{ for $1\leq i \leq \frac{n}{2}-1$}
\end{align} 

In order to simplify the analysis, for every $j=0, \ldots \frac{n}{2}-1$, we write 
\begin{align*}
C_{\frac{n}{2}-j}= \alpha \frac{P_j(\alpha)}{Q_j(\alpha)}
\end{align*}
where $P_j(\alpha)$ and $Q_j(\alpha)$ are polynomials in $\alpha$. Using the fact that $C_{\frac{n}{2}}=2\alpha$, we find that $P_0(\alpha)=2$ and $Q_0(\alpha)=1$, while 
using Equation~\ref{eq:C2} we find that for every $j=1, \ldots \frac{n}{2}-1$
\begin{align*}
P_j(\alpha)&= Q_{j-1}(\alpha)
\\
Q_j(\alpha)&= Q_{j-1}(\alpha) - \alpha^2 P_{j-1}(\alpha)
\end{align*}
which gives the following second order recurrence relation for $Q_j(\alpha)$
\[
Q_j(\alpha)= Q_{j-1} - \alpha^2 Q_{j-2}\;\qquad \text{ with $Q_0(\alpha)=1$ and $Q_1(\alpha)=1-2\alpha^2$}
\]
Consider the  characteristic equation of the second order recurrence relation for $Q_j$, i.e., 
\[
x^2 -x +\alpha^2=0
\]
whose  solutions are $x_1= \frac{1 + \sqrt{1-4\alpha^2}}{2}$ and $x_2= \frac{1 - \sqrt{1-4\alpha^2}}{2}$.
Then, we know that $Q_j(\alpha)$ satisfies the following equation
\[
Q_j(\alpha)= A x_1^j + B x_2^j
\]
where $A$ and $B$ can be computed using the initial conditions $Q_0(\alpha)=1$ and $Q_1(\alpha)=1-2\alpha^2$. Specifically, we obtain that $A=x_1$ and $B=x_2$. Overall, we have that for every $j=0, \ldots \frac{n}{2}-1$
\[
Q_j(\alpha)= x_1^{j+1} +  x_2^{j+1}
\]
and 
\begin{align*}
C_1 &= C_{\frac{n}{2}- \left(\frac{n}{2}-1\right)}= \alpha 
\frac{P_{\frac{n}{2}-1}(\alpha)}{Q_{\frac{n}{2}-1}(\alpha)}= \alpha 
\frac{Q_{\frac{n}{2}-2}(\alpha)}{Q_{\frac{n}{2}-1}(\alpha)} = \alpha \frac{x_1^{\frac{n}{2}-1} +  x_2^{\frac{n}{2}-1}}{x_1^{\frac{n}{2}} +  x_2^{\frac{n}{2}}}
\\
&= \frac{\alpha}{x_1}\; \frac{1  + \left( \frac{x_2}{x_1} \right)^{\frac{n}{2}-1}   }{ 1  + \left( \frac{x_2}{x_1} \right)^{\frac{n}{2}} }
\end{align*}

\bibliographystyle{abbrv}
\bibliography{RW-epidemics}

\end{document}